\documentclass[11pt]{amsart}
\usepackage[a4paper,margin=2.6cm]{geometry}
\usepackage{amsmath,amssymb,amsthm}
\usepackage{xcolor}
\usepackage[hidelinks]{hyperref}
\newcommand{\rev}[1]{{#1}}
\usepackage[T1]{fontenc}
\usepackage{lmodern}
\hypersetup{pdftitle={Boundary control of longitudinal vibrations in composite beams with localized shear connectors},pdfauthor={Kaïs Ammari and Mourad Sini},bookmarksnumbered=true}
\allowdisplaybreaks[1]
\newtheorem{theorem}{Theorem}[section]
\newtheorem{proposition}[theorem]{Proposition}
\newtheorem{lemma}[theorem]{Lemma}

\theoremstyle{definition}

\newtheorem{example}[theorem]{Example}
\theoremstyle{remark}

\newtheorem*{remark*}{Remark}

\numberwithin{equation}{section}

\newcommand{\R}{\mathbb{R}}

\title[Boundary control of connected beams]{Boundary control of longitudinal vibrations in composite beams with localized shear connectors}
\author{Ka\"is Ammari}
\address{LR Analysis and Control of PDEs, LR22ES03, Department of Mathematics, Faculty of Sciences of Monastir, University of Monastir, Tunisia}
\email{kais.ammari@fsm.rnu.tn}
\author{Mourad Sini}
\address{RICAM, Austrian Academy of Sciences, Altenbergerstrasse 69, 4040 Linz, Austria}
\email{mourad.sini@oeaw.ac.at}
\date{}
\newcommand{\X}{\mathcal X}
\newcommand{\U}{\mathcal U}
\newcommand{\Tad}{\Theta_{\mathrm{ad}}}
\newcommand{\mass}{\mathsf M}
\newcommand{\form}{\mathfrak a}
\newcommand{\op}{\mathsf A}
\newcommand{\refn}{\mathrm b}
\begin{document}
\begin{abstract}
\noindent We are concerned with the problem of controlling the longitudinal vibrations of a composite beam formed
by two beams connected through a finite number of localized shear connectors. Our aim is to understand
how the positions and the stiffnesses of the connectors affect the possibility, and the cost, of driving
the system to rest by means of axial forces applied at one end. We first justify the point-connector
model by approximating the shear stiffness of each connector by a density supported in an interval of
length $2\varepsilon$ and of prescribed integral. We derive an estimate of order $\varepsilon^{1/2}$ for
the difference of the elastic forms and of the corresponding spatial resolvents, and establish strong
convergence of the free evolutions in the energy space. Let $L$ be the length of the beams and $c_1,c_2$
their longitudinal wave speeds. In the case of two independent boundary forces, we prove exact
controllability in the energy space if and only if $T\geq 2L/\min(c_1,c_2)$. In contrast, with only one
boundary force, exact controllability in the same space fails for every finite time, owing to the
compactness of the transfer to the uncontrolled beam. \rev{From rest, the reachable set of the second member
is trivial up to the first arrival at an active connector and is infinite-dimensional and nonclosed
thereafter.} We also derive an observability estimate,
including at the critical time, which is uniform with respect to the positions and stiffnesses of the
connectors in a compact admissible family. As a consequence, under a constraint on the total stiffness
and separation conditions on the connectors, we prove the existence of an optimal configuration
minimizing the worst control cost in the presence of \rev{simultaneous partial or total loss of stiffness in up to
a prescribed number of connectors}. Finally, for a prescribed finite frequency band, we characterize the controllability by
one force and the corresponding cost in terms of a finite-dimensional Gram matrix.
\end{abstract}
\maketitle
\setlength{\parindent}{0pt}
\setlength{\parskip}{3pt plus 1pt minus 1pt}
\section{Introduction}
\subsection{Motivation and the mechanical model}
\noindent Consider a composite beam formed by the connection of a steel beam and a reinforced concrete slab. One
of the relevant quantities in describing the integrity of such a structure is the shearing stiffness
of the connection, which accounts for the resistance to the relative longitudinal displacement of the
two members. The role of the connection in the dynamical behaviour of composite beams has been
investigated both experimentally and mathematically. We refer, in particular, to Biscontin, Morassi
and Wendel~\cite{BMW} for experiments with different densities of connectors and to Morassi and
Rocchetto~\cite{MR} for the influence of damage on the vibration modes. In~\cite{MNS}, the determination
of the shearing stiffness from dynamical measurements at one end of the composite beam was addressed.

\noindent In the present work, we look at this system from the point of view of control. More precisely, given a prescribed number of connectors, how should we distribute them, and assign their stiffnesses, in such a way that the longitudinal vibrations of both members can be suppressed by forces applied at one end, at a cost which remains as small as possible \rev{when several connectors can be damaged simultaneously}? There are
two issues to be distinguished here: the controllability of the system, for a given configuration of connectors, and the choice of a configuration in view of the cost of the control. We shall also compare the cases where both members and where only one of them can be acted upon at the accessible end.

\noindent We restrict ourselves to longitudinal vibrations and use the mechanical model considered
in~\cite[eq.~(1.4)]{MNS}, in which the bending motions are neglected. In contrast to a continuously
distributed connection, we suppose that the shearing stiffness is concentrated in small neighbourhoods
of a finite number of points of the beam. The point model used below is justified in
Section~\ref{sec:forms} as the limit of such localized distributions.

\noindent Let $L>0$ denote the length of the beam and $u_i(x,t)$ the axial displacement of the $i$-th member,
$i=1,2$, at $(x,t)\in(0,L)\times(0,T)$. Denote by $\rho_i>0$ the mass per unit length and by
$a_i=E_i A_i>0$ the axial rigidity (the product of Young's modulus and the cross-sectional area). For
simplicity, these coefficients are taken to be constant. We set
\[
 s_u=u_1-u_2,\qquad b=(1,-1)^{\mathsf t},\qquad
 \mass=\operatorname{diag}(\rho_1,\rho_2),\qquad D=\operatorname{diag}(a_1,a_2).
\]
For a nonnegative stiffness density $k(x)$, the elastic energy associated with the connection is
$\frac12\int_0^L k(x)|s_u(x)|^2\,dx$. If the stiffness is concentrated near a point $\xi_j$ and its
integral is $k_j$, the limiting energy is $\frac12 k_j|s_u(\xi_j)|^2$. Thus the density $k(x)$ is a
shearing stiffness per unit length, whereas $k_j$ is the stiffness of a single connector, measured
as force per unit slip.

\noindent Let $M\geq1$ be the number of nominal connectors. Given $0<\xi_1<\cdots<\xi_M<L$ and $k_j\geq0$, we set
$\theta=(\xi,k)$ and $\nu_\theta=\sum_{j=1}^M k_j\delta_{\xi_j}$. A vanishing stiffness represents a
completely failed connector. We consider the following initial-boundary value problem:
\begin{equation}\label{eq:model}
\begin{aligned}
 \mass u_{tt}-(D u_x)_x+b\,s_u\,\nu_\theta&=0,\\
 u(L,t)=0,\qquad -D u_x(0,t)&=p(t),\\
 (u,u_t)|_{t=0}&=z .
\end{aligned}
\end{equation}
Here, the two members are fixed at $x=L$ and $p_i=-a_i u_{i,x}(0,\cdot)$ is the axial force applied to
the $i$-th member at $x=0$. We measure the cost of a control by the norm $\|p\|_{L^2(0,T;\mathbb R^2)}$. Our main purpose is to choose the nominal positions and stiffnesses of the connectors so as to minimize the worst cost of driving the system to rest at time $T$. To take into account the possible deterioration of the connection, \rev{we allow up to a prescribed number of connectors to lose part or all of their stiffness at the same time}. The configuration is chosen before the damage occurs, whereas the control is chosen
for the actual (known) damaged configuration and the initial state. The stiffnesses are assumed to
remain constant during each control experiment.
\subsection{Functional setting and statement of the main results}
\noindent We first introduce some notations and the functional framework of the problem. We work with real spaces, using their complexifications when dealing with spectral arguments. We set
\[
 \begin{aligned}
 H&=L^2(0,L)^2, & (v,w)_H&=\int_0^L(\mass v)\cdot w\,dx,\\
 V&=\{v\in H^1(0,L)^2:v(L)=0\}, & \X&=V\times H .
 \end{aligned}
\]
To compare different configurations, we use on $\X$ the reference norm
$\|z\|_{\refn}^2=\form_{\refn}(z_0,z_0)+\|z_1\|_H^2$, where
$\form_{\refn}(u,v)=\int_0^L Du'\cdot v'\,dx$ and $\|u\|_V^2=\form_{\refn}(u,u)$. The subscript $\refn$
refers to the uncoupled system. Given a finite nonnegative stiffness measure $\nu$, we introduce the
elastic form and the associated energy scalar product by
\begin{equation}\label{eq:form}
 \form_\nu(u,v)=\form_{\refn}(u,v)+\int s_us_v\,d\nu,\qquad
 (z,y)_\nu=\form_\nu(z_0,y_0)+(z_1,y_1)_H .
\end{equation}
We denote by $\|\cdot\|_\nu$ the norm induced by the latter scalar product, by $\op_\nu$ the spatial
operator in $H$ associated with $\form_\nu$, and by $S_\nu(t)$ the free wave group on $\X$. As a subscript,
$\theta$ stands for $\nu_\theta$ and $\refn$ for the zero measure. To approximate the point connectors,
we take a nonnegative function $\eta\in C_c^\infty(-1,1)$ of integral one and set
\begin{equation}\label{eq:foot}
 d\nu_{\theta,\varepsilon}(x)=\sum_{j=1}^M
 k_j\varepsilon^{-1}\eta((x-\xi_j)/\varepsilon)\,dx .
\end{equation}
In what follows, $(\theta,\varepsilon)$, as a subscript, refers to the measure $\nu_{\theta,\varepsilon}$. We use $C$ to denote a generic positive constant whose value may change from one estimate to another and which is independent of the parameters under consideration, unless
otherwise stated. We start with the result which justifies the use of point connectors.

\begin{theorem}[Point-connector limit]\label{thm:limit}
Let $k_j\geq0$, $j=1,\ldots,M$. The free problem associated with \eqref{eq:model} generates a group of
isometries on $(\X,\|\cdot\|_\theta)$, with conserved energy $\frac12\|(u,u_t)\|_\theta^2$. The domain
of the spatial operator, including the transmission conditions at the connectors, is given in
Proposition~\ref{prop:domain}. Assume that $[\xi_j-\varepsilon,\xi_j+\varepsilon]\subset(0,L)$ for all $j$.
Then, for every $\lambda>0$, we have
\[
 \|\form_{\theta,\varepsilon}-\form_\theta\|_{V\times V}
  +\|(\op_{\theta,\varepsilon}+\lambda)^{-1}
            -(\op_\theta+\lambda)^{-1}\|_{H\to V}
 \le C_\lambda\sqrt{\varepsilon}.
\]
Here the norm on the space of bilinear forms is the supremum over the two unit balls in $V$.
Moreover, for every $z\in\X$ and every $R>0$,
\[
 \sup_{|t|\le R}\|S_{\theta,\varepsilon}(t)z-S_\theta(t)z\|_{\refn}\longrightarrow0
 \quad\text{as }\varepsilon\longrightarrow0.
\]
The form and resolvent estimates above are uniform with respect to the configurations with bounded
total stiffness.
\end{theorem}

\noindent Having justified the mechanical model, we address the controllability problem for a fixed configuration.
\begin{theorem}[Sharp time for boundary controllability]\label{thm:control}
Set $c_i=\sqrt{a_i/\rho_i}$, $i=1,2$, and $T_*=2L/\min_i c_i$. For every finite configuration of
connectors of nonnegative stiffnesses, every $T>0$, every $p\in\U_T:=L^2(0,T;\R^2)$ and every $z\in\X$,
problem \eqref{eq:model} admits a unique energy solution. In addition, the following properties hold.
\begin{enumerate}
 \item \noindent With two independent boundary forces, any two states in $\X$ can be joined in time $T$ if and
 only if $T\geq T_*$. For such times, the unique control of minimum norm is given by \eqref{eq:HUM}.
 \item \noindent \rev{Suppose that $p_2=0$. For the system starting from rest, denote by
 $K_{\theta,T}p_1=\pi_2R_{\theta,T}(p_1,0)$ the terminal displacement and velocity of the second
 member, where $R_{\theta,T}$ is the endpoint map and $\pi_2$ is the projection onto that member.
 Let
 \[
 t_a=\min_{j:\,k_j>0}\frac{\xi_j}{c_1},
 \]
 with the convention that $ \min\varnothing=+\infty$. The operator $K_{\theta,T}$ is compact for every finite $T>0$. Its range is $\{0\}$ for
 $T\leq t_a$ and is infinite-dimensional and nonclosed for $T>t_a$. Consequently, exact
 controllability in $\X$ by one boundary force fails for every finite time.}
\end{enumerate}
\end{theorem}

\noindent We now formulate the optimization problem for the positions and stiffnesses of the connectors. We fix
their number $M$, a minimum separation and boundary clearance $d>0$, the total nominal stiffness
$K>0$ and bounds $0\leq k_-\leq k_+<\infty$ on the individual nominal stiffnesses. We assume that the
following set is nonempty:
\begin{equation}\label{eq:admiss}
\begin{split}
 \Tad=\{(\xi,k):\ &d\le \xi_1,\quad \xi_M\le L-d,\quad
 \xi_{j+1}-\xi_j\ge d\ (j<M),\\
 & k_-\le k_j\le k_+,\quad \textstyle\sum_j k_j=K\}.
\end{split}
\end{equation}

\noindent \rev{Fix an integer $1\leq r\leq M$ and a minimum retained stiffness fraction
$\gamma_0\in[0,1]$. The $j$-th connector retains stiffness $\gamma_j k_j$. For a connector of
positive nominal stiffness, $\gamma_j=1$ describes the intact state, $0<\gamma_j<1$ partial
damage, and $\gamma_j=0$ complete failure. We allow damage to at most $r$ connectors and set
\[
 \Gamma_r=\bigl\{\gamma\in[\gamma_0,1]^M:
           \#\{j:\gamma_j<1\}\leq r\bigr\},\qquad
 \theta^\gamma=\bigl(\xi,(\gamma_j k_j)_{j=1}^M\bigr).
\]
The vector $(1,\ldots,1)$ represents the intact configuration and belongs to $\Gamma_r$ for every $r$. The total stiffness
constraint applies to the nominal stiffnesses, so that the total stiffness after damage is
at most $K$.} For $T\geq T_*$, we define
\begin{equation}\label{eq:design}
\begin{aligned}
 m_T(\theta,z)&=\min\{\|p\|_{\U_T}:(u,u_t)(0)=z,\ (u,u_t)(T)=0\},\\
 \mathcal C_T(\theta)&=\sup_{\|z\|_{\refn}=1}m_T(\theta,z),\qquad
 \rev{\mathcal J_{T,r}(\theta)=\sup_{\gamma\in\Gamma_r}\mathcal C_T(\theta^\gamma).}
\end{aligned}
\end{equation}
The quantity $m_T(\theta,z)$ is the least force norm needed to drive the initial state $z$ to rest,
$\mathcal C_T(\theta)$ is the corresponding cost for the worst unit initial state, and
\rev{$\mathcal J_{T,r}(\theta)$ is the worst cost over all the damage scenarios in $\Gamma_r$.} We have the following result.

\begin{theorem}[Existence of a robust configuration]\label{thm:design}
Let $T\geq T_*$ \rev{and $1\leq r\leq M$}. The observability constant associated with the two boundary forces admits a strictly
positive lower bound, uniform with respect to all the intact and damaged configurations entering
\eqref{eq:design}. Moreover, $\mathcal C_T$ is lower semicontinuous on this family and is given by the
inverse square root of the optimal observability constant with the dual reference normalization in
\eqref{eq:best}. In particular, there exists $\theta_*\in\Tad$ such that
\[
 \rev{\mathcal J_{T,r}(\theta_*)=\min_{\theta\in\Tad}\mathcal J_{T,r}(\theta)<\infty.}
\]
\rev{For $r=1$ this includes damage to any single connector. If $\gamma_0=0$,
any subset of at most $r$ connectors may fail completely, including all connectors when $r=M$.}
\end{theorem}

\subsection{Interpretation of the results in view of the original problem}\label{sec:interpretation}

\noindent
\rev{A shear connector resists the relative longitudinal displacement of the
two members. Theorem~\ref{thm:limit} justifies the point-connector model for
the free dynamics when the width of each connector decreases while its
integrated stiffness $k_j$ remains fixed. In this limit, the displacements
remain continuous across $\xi_j$, and the axial forces have opposite jumps
proportional to the local slip $s_u(\xi_j,t)$, as expressed in
Proposition~\ref{prop:domain}. Thus $\xi_j$ determines where force can be
transferred between the members, while $k_j$ determines the magnitude
$k_j|s_u(\xi_j,t)|$ of the force transmitted for a given slip. The design
problem is to choose these positions and distribute the available stiffness
$\sum_{j=1}^M k_j=K$ among them.}

\noindent
\rev{Suppose first that an independent axial force can be applied to each
member at $x=0$. By Theorem~\ref{thm:control}, the shortest time in which
an arbitrary initial vibration can be brought to rest is
$T_*=2L/\min(c_1,c_2)$, the time required for a longitudinal wave in the
slower member to travel to the fixed end and return. Moving a connector
closer to the loaded end or increasing its stiffness cannot shorten this
time. Even if all connectors fail, both members can still be brought to
rest for $T\geq T_*$, since each has its own actuator. For a prescribed
time $T\geq T_*$, the choice of connectors therefore concerns the amount of
boundary force required to suppress the vibration. In the present model,
this effort is measured by $\|p\|_{\U_T}$.}

\noindent
\rev{With a single force applied to the first member, the connectors provide
the only means of transmitting motion to the second. Starting from rest,
the second member remains at rest up to the first-arrival time
$t_a=\min_{j:\,k_j>0}\xi_j/c_1$. Moving the first active connector closer
to the loaded end advances the earliest possible transfer. If that
connector fails completely, the next connector of positive stiffness
determines the new arrival time. If none remains, the two members are
disconnected and $t_a=+\infty$. For $T>t_a$, the reachable states of the
second member form an infinite-dimensional, nonclosed subspace of its
energy space. Thus the actuator can generate infinitely many independent
vibration profiles in that member, although it cannot prescribe an
arbitrary terminal displacement and velocity there.}

\noindent
\rev{A natural objective with one actuator is to suppress vibrations in a
prescribed frequency band. Example~\ref{ex:node} illustrates why the
connector positions matter for this purpose. If the only active connector
lies at a node of a mode supported on the second member, the slip at that
connector vanishes for this mode. Increasing its stiffness then has no
effect on the mode, which remains invisible at the controlled end.
For a prescribed band, Proposition~\ref{prop:band} gives the criterion
for eliminating its modal components and quantifies the required control
effort through the associated Gram matrix.}

\noindent
\rev{Returning to the case of two boundary forces, Theorem~\ref{thm:design}
gives the design criterion a precise interpretation in the presence of
damage. The nominal configuration $\theta\in\Tad$ is chosen before damage
occurs, subject to the stiffness budget and separation conditions
in~\eqref{eq:admiss}. For each retained stiffness vector $\gamma\in\Gamma_r$,
the force can then be adapted to the actual configuration $\theta^\gamma$
and to the initial vibration $z$. Since
$\mathcal C_T(\theta^\gamma)=\beta_T(\theta^\gamma)^{-1/2}$, minimizing
the worst control cost in~\eqref{eq:design} is equivalent to solving}
\begin{equation}\label{eq:robust-beta}
 \rev{\max_{\theta\in\Tad}\ \inf_{\gamma\in\Gamma_r}
       \beta_T(\theta^\gamma).}
\end{equation}
\rev{For an optimal nominal configuration $\theta_*$, every admissible
damage scenario and every initial state $z\in\X$ admit a force $p$ that
brings both members to rest at time $T$ and satisfies}
\begin{equation}\label{eq:robust-force}
 \rev{\|p\|_{\U_T}\leq \mathcal J_{T,r}(\theta_*)\,\|z\|_{\refn}.}
\end{equation}
\rev{The coefficient $\mathcal J_{T,r}(\theta_*)$ is the smallest such
worst-case bound among the admissible designs. In particular, if $r=M$
and $\gamma_0=0$, complete disconnection is an admissible scenario for
every design. The cost of stopping the two disconnected members then
gives a lower bound independent of the nominal connector positions and
stiffnesses, as calculated in Example~\ref{ex:freecost}. The theorem
establishes the existence of an optimal arrangement for active vibration
control. Determining its positions and stiffnesses for a particular beam
remains a quantitative design problem.}

\subsection{Related works and organization of the paper}
\noindent Let us relate our results to some of the works on connected beams and on the control of systems with
several wave speeds. Belishev, Blagovestchenskii and Ivanov~\cite{BBI} studied boundary controllability
for a two-velocity system and used it in the analysis of inverse problems. In the case of connected
beams with distributed shearing stiffness, Belishev and Ivanov~\cite{BI} recovered this coefficient
from the response at the loaded end together with a measurement at the fixed end. In~\cite{MNS}, the
reconstruction was established from the frequency response at one end. The model including bending
motions, and hence a further coefficient of the connection, was investigated in~\cite{MNSS}, using
boundary and interior data. We also refer to Li and Rao~\cite{LiRao} for exact controllability of a
coupled wave system with Neumann controls and for the lack of exact controllability in the common
energy space when the number of controls is reduced.

\noindent Here the main points are the treatment of the concentrated stiffnesses, the inclusion of the critical
control time and the uniformity of the observability estimate with respect to the configuration of
the connectors. To obtain this estimate, we propagate the observed boundary velocities along the
characteristics, taking into account the jumps at the connectors. At the middle of the critical time
interval, the resulting characteristic region covers the entire spatial interval. A compactness and
uniqueness argument, allowing the connector positions to vary and the stiffnesses to vanish, then
yields the uniform observability estimate. This, together with the lower semicontinuity of the cost,
is the key ingredient in the proof of Theorem~\ref{thm:design}.
\bigskip

\noindent The rest of the paper is organized as follows. In Section~\ref{sec:forms}, we justify the point model
and prove Theorem~\ref{thm:limit}. Section~\ref{sec:forces} is devoted to the well-posedness of the controlled system and to a compactness property of its endpoint map. In Sections~\ref{sec:two}
and~\ref{sec:one}, we prove the assertions of Theorem~\ref{thm:control}, for two and one boundary
forces, respectively. Finally, in Section~\ref{sec:design}, we prove Theorem~\ref{thm:design} and
address the corresponding design problem for the control of a finite frequency band by a single
force.
\section{From distributed stiffness to point connectors}\label{sec:forms}
\subsection{Closed form and axial-force transmission}
\noindent We start by recalling some estimates for the slip, which will be used throughout the paper. Set
$\alpha=(a_1^{-1}+a_2^{-1})^{1/2}$. By the Cauchy--Schwarz inequality applied to $s'_u$ and then to its
integral between two points, we have
\begin{equation}\label{eq:holder}
 |s_u(x)-s_u(y)|\le\alpha |x-y|^{1/2}\|u\|_V,\qquad
 \|s_u\|_\infty\le\alpha\sqrt L\,\|u\|_V .
\end{equation}
The second estimate follows from the first one by taking $y=L$. By Poincar\'e's inequality, $V$ is
continuously embedded in $H$. It follows from \eqref{eq:holder} that, for $\nu([0,L])\leq K$,
\begin{equation}\label{eq:uniform}
 \|u\|_V^2\le\form_\nu(u,u)\le (1+K\alpha^2 L)\|u\|_V^2 .
\end{equation}
In particular, $\form_\nu$ is a densely defined, symmetric and closed form on $H$. Denote by $\op_\nu$
the positive self-adjoint operator associated with this form, i.e., $u\in D(\op_\nu)$ if and only if
$u\in V$ and there exists $g\in H$ such that $\form_\nu(u,v)=(g,v)_H$ for every $v\in V$. In this case,
$\op_\nu u=g$. We have $D(\op_\nu^{1/2})=V$. Moreover, $\op_\nu$ is strictly positive by Poincar\'e's
inequality and its resolvent is compact by the compact embedding of $V$ in $H$. We give next the
domain of the operator in the case of point connectors.

\begin{proposition}\label{prop:domain}
Set $\xi_0=0$, $\xi_{M+1}=L$ and denote by $[w]_{\xi}=w(\xi^+)-w(\xi^-)$ the jump of $w$ at $\xi$.
The domain $D(\op_\theta)$ consists of the functions $u\in V$ such that
$u|_{(\xi_j,\xi_{j+1})}\in H^2((\xi_j,\xi_{j+1}))^2$, $j=0,\ldots,M$, and
\[
 Du'(0)=0,\qquad [Du']_{\xi_j}=k_j b\,s_u(\xi_j).
\]
On each subinterval $\op_\theta u=-\mass^{-1}Du''$. In particular the sum of the two axial forces is continuous at a connector.
\end{proposition}
\begin{proof}
Let $u\in D(\op_\theta)$ and set $g=\op_\theta u$. By taking test functions supported in the open
subintervals, we obtain $-Du''=\mass g$ there. Since $g\in H$ and $u\in H^1$, the restriction of $u$
to each subinterval belongs to $H^2$. Integrating by parts in the variational identity, we obtain
\[
 -Du'(0)\cdot v(0)+\sum_j\{k_j b\,s_u(\xi_j)-[Du']_{\xi_j}\}\cdot v(\xi_j)=0.
\]
Since the values of $v$ at $0,\xi_1,\ldots,\xi_M$ can be prescribed independently (take, for instance,
continuous piecewise affine functions in $V$), we obtain the stated conditions. Conversely, these
conditions imply the variational identity defining $\op_\theta$, by integration by parts. The last
assertion follows on summing the two components of the jump, since $b_1+b_2=0$.
\end{proof}

\noindent We now consider the evolution associated with $\op_\nu$. Introduce the operator
$G_\nu(u,v)=(v,-\op_\nu u)$ with domain $D(G_\nu)=D(\op_\nu)\times V$. From the defining identity of the
form it follows that $G_\nu$ is skew-symmetric with respect to $(\cdot,\cdot)_\nu$. Moreover, for $h>0$,
the operators $hI\pm G_\nu$ are surjective as eliminating the second component reduces the equation to
a variational problem with coercive form $\form_\nu+h^2(\cdot,\cdot)_H$ (the resolvent formulas are
given below). The first component of the solution belongs to $D(\op_\nu)$ by the variational identity.
Hence $G_\nu$ is skew-adjoint and generates a group $S_\nu$ of isometries. Equivalently, by spectral
calculus, the displacement is given by
\[
 u(t)=\cos(t\sqrt{\op_\nu})u(0)
       +\op_\nu^{-1/2}\sin(t\sqrt{\op_\nu})u_t(0).
\]
For initial data in $\X$, the displacement satisfies $u\in C(\R;V)\cap C^1(\R;H)$ and
$u_{tt}\in C(\R;V')$, where $H$ is identified with its dual via $(\cdot,\cdot)_H$. In the point-connector
model the conserved energy takes the form
\[
 E_\theta(t)=\tfrac12\sum_{i=1}^2\int_0^L(\rho_i|u_{i,t}|^2+a_i|u_i'|^2)\,dx
                 +\tfrac12\sum_j k_j|s_u(\xi_j)|^2 .
\]

\subsection{Approximation by distributed stiffnesses in a common energy space}
\noindent \rev{For scalar one-dimensional operators, Brasche, Figari and Teta~\cite{BFT} proved that weak
convergence of finite potential measures gives norm-resolvent convergence of the associated
Schr\"odinger operators. Bertini, Noja and Posilicano~\cite{BNP} studied wave equations with point
interactions in finite energy spaces. A related three-dimensional wave--oscillator model is treated
in~\cite{CGMS}. For the present two-member system on a bounded interval, the embedding
$V\hookrightarrow C([0,L])^2$ gives a direct estimate of order $\sqrt{\varepsilon}$ for the
difference of the elastic forms.} Set
$m_{1/2}=\int |r|^{1/2}\eta(r)\,dr$. Adding and subtracting $s_u(\xi_j)s_v(x)$ under the integral and
using \eqref{eq:holder}, we obtain
\begin{equation}\label{eq:formerror}
 \left|\form_{\theta,\varepsilon}(u,v)-\form_\theta(u,v)\right|
 \le 2\alpha^2\sqrt L\,m_{1/2}\Big(\sum_j k_j\Big)
       \sqrt{\varepsilon}\,\|u\|_V\|v\|_V.
\end{equation}
Notice that the normalization by $\varepsilon^{-1}$ in \eqref{eq:foot} is essential here: it keeps
the integrated stiffness of the $j$-th connector equal to $k_j$, independently of $\varepsilon$.

\noindent For the design problem, we shall also need to compare two point configurations. Let $\theta$ and
$\tilde\theta$ have $M$ connectors each and total stiffness bounded by $K$. Setting
$h_\xi=\max_j|\xi_j-\tilde\xi_j|$, we obtain, in a similar manner,
\begin{equation}\label{eq:moving}
 \|\form_\theta-\form_{\tilde\theta}\|_{V\times V}
 \le \alpha^2 L\sum_j|k_j-\tilde k_j|
                 +2\alpha^2\sqrt L\,K\sqrt{h_\xi}.
\end{equation}
The first term accounts for the variation of the stiffnesses at fixed positions, while the second
one follows by applying \eqref{eq:holder} to the positions $\xi_j$ and $\tilde\xi_j$. Note that these
estimates remain valid when some of the stiffnesses vanish.

\begin{proof}
\rev{We prove the resolvent and wave convergence in Theorem~\ref{thm:limit}.}
Denote by $B_\nu:V\to V'$ the operator associated with the form and by $J:H\to V'$ the canonical
embedding, i.e., $\langle Jg,v\rangle=(g,v)_H$. Given $F\in V'$, set
$w_\nu=(B_\nu+\lambda J|_V)^{-1}F$. By coercivity, $\|w_\nu\|_V\leq\|F\|_{V'}$. Subtracting the
variational identities corresponding to $\nu$ and $\nu'$ and testing with $w_{\nu'}-w_\nu$, we obtain
\[
 \|w_{\nu'}-w_\nu\|_V\le
       \|\form_{\nu'}-\form_\nu\|_{V\times V}\|F\|_{V'}.
\]
Taking $F=Jg$, $g\in H$, and using \eqref{eq:formerror} gives the stated estimate for the spatial
resolvents, in particular their convergence in the operator norm on $H$. To pass to the wave groups,
we work on the fixed space $(\X,\|\cdot\|_{\refn})$. For $h>0$ and $(f,g)\in\X$, we have
\[
 (hI-G_\nu)^{-1}(f,g)=(w,hw-f),\quad
 w=(B_\nu+h^2J|_V)^{-1}J(g+hf).
\]
Hence the resolvents of the generators converge in the operator norm on $(\X,\|\cdot\|_{\refn})$.
For $-G_\nu$ the analogous formula is $w=(B_\nu+h^2J|_V)^{-1}J(hf-g)$, with second component $f-hw$.
The groups and their inverses are uniformly bounded in the reference norm by \eqref{eq:uniform}.
Applying the Trotter--Kato approximation theorem~\cite[Ch.~III, Thm.~4.8]{EN} to both signs, we deduce
the asserted strong convergence of the groups, uniformly on compact time intervals. This completes
the proof of Theorem~\ref{thm:limit}. In addition, the same arguments, using \eqref{eq:moving} instead
of \eqref{eq:formerror}, show that $S_{\theta_n}(t)z\to S_\theta(t)z$ uniformly on compact time
intervals for each $z\in\X$ whenever $\theta_n\to\theta$.
\end{proof}
\medskip\noindent The above estimate is of order $\varepsilon^{1/2}$ at the level of the forms and resolvents. For the free evolutions, we have obtained strong convergence for each fixed datum, without asserting a rate in the energy norm. To justify the approximation of the optimal controls themselves, one would need, in addition, to study the boundary observations and their observability estimates for the distributed
stiffness model. This question is beyond the scope of the present work.

\section{Boundary forces and energy solutions}\label{sec:forces}
\noindent In this section, we establish the well-posedness of \eqref{eq:model} with $L^2$ boundary forces. The solution is defined by transposition using the velocity at the accessible end of solutions of the
homogeneous adjoint problem. To make sense of this formulation in the energy space, we first derive
an estimate for the velocity trace.
\subsection{Direct trace estimates}
\noindent For a homogeneous solution $\psi$, set $\Psi=(\psi,\psi_t)$. We introduce the observation operator,
initially for regular initial data $y\in\X$, by
\[
 (O_{\theta,T}y)(t)=\psi_t(0,t),\qquad \Psi(t)=S_\theta(t)y .
\]
The following lemma allows us to extend this definition to all $y\in\X$, with traces understood in
$L^2$ in time.
\begin{lemma}\label{lem:trace}
Let $T>0$. The operator $O_{\theta,T}:\X\to\U_T$ is bounded. Its norm is bounded uniformly with
respect to the configurations with total stiffness at most $K$ and $\xi_1\geq d>0$. In addition,
for each connector of a fixed configuration, the map $y\mapsto b\cdot\psi_t(\xi_j,\cdot)$ is bounded
from $\X$ into $L^2(0,T)$.
\end{lemma}
\begin{proof}
We first prove the estimates for initial data given by finite linear combinations of eigenfunctions
of $\op_\theta$, which are dense in $\X$. For such solutions, both $\psi$ and $\psi_t$ are continuous
across each connector. On an interval $(A,B)$ containing no connector, the component $w=\psi_i$
satisfies $\rho_i w_{tt}=a_i w_{xx}$. Set $F_i=\rho_i w_t^2+a_i w_x^2$. Multiplying the equation by
$q(x)w_x$, where $q$ is a smooth function, and integrating in space and time, we obtain
\begin{equation}\label{eq:mult}
 \tfrac12\int_0^\tau [qF_i]_A^B\,dt
 =\left[\rho_i\int_A^B q w_t w_x\,dx\right]_0^\tau
       +\tfrac12\int_0^\tau\!\int_A^B q'F_i\,dx\,dt .
\end{equation}
To estimate the trace at $0$, take $A=0$, $B=d/2$ and choose $q$ such that $q(0)=1$ and $q(B)=0$
(for a fixed configuration any $0<d\leq\xi_1$ suffices). By conservation of energy, the absolute
value of the right-hand side is bounded by
$C(\|q\|_\infty+\tau\|q'\|_\infty)E_\theta(0)$. To estimate the trace at $\xi_j$, we take the interval
immediately to its left and choose $q(A)=0$, $q(\xi_j)=1$, using the derivative on that side. Since
$F_i\geq\rho_i w_t^2$, the estimates follow. By density, we extend the velocity traces to all
initial data in $\X$. Moreover, for a sequence converging in the energy space, the displacements at
any fixed point converge uniformly in time by the spatial trace inequality. Hence the velocity trace
so defined is the distributional time derivative of the displacement trace. For $0\leq\tau\leq T$,
the constants at the accessible end can be chosen uniformly under the stated clearance condition.
\end{proof}

\subsection{Transposition and the endpoint map}
\noindent For $0\leq\tau\leq T$, denote by $R_{\theta,\tau}p\in\X$ the state at time $\tau$ corresponding to
zero initial data and the restriction of $p$ to $(0,\tau)$. We define it by the transposition identity
\begin{equation}\label{eq:reach}
 (R_{\theta,\tau}p,h)_\theta
 =\int_0^\tau p(t)\cdot\psi_t^h(0,t)\,dt,\qquad
 \Psi^h(t)=S_\theta(t-\tau)h,\quad h\in\X .
\end{equation}
The sign in this identity follows from integration by parts in space: the boundary contribution
at $x=0$ to the derivative of the energy is $-Du_x(0,t)\cdot u_t(0,t)=p(t)\cdot u_t(0,t)$. Since
$S_\theta(-\tau)$ is an isometry for the energy norm, Lemma~\ref{lem:trace} also applies to the
solutions with prescribed final data. By the Riesz representation theorem, there exists a unique
state satisfying \eqref{eq:reach}. Furthermore, using \eqref{eq:uniform}, we have
\begin{equation}\label{eq:energybd}
 \sup_{0\le\tau\le T}\|S_\theta(\tau)z+R_{\theta,\tau}p\|_{\refn}
       \le C(T,d,K)(\|z\|_{\refn}+\|p\|_{\U_T}) .
\end{equation}
Here and in what follows the constants may depend also on the material coefficients and on $L$.
Splitting the integral in \eqref{eq:reach} at time $s$ gives
\[
 R_{\theta,s+h}p=S_\theta(h)R_{\theta,s}p+
 R_{\theta,h}\big(p(s+\cdot)\big),\qquad 0\le s\le s+h\le T.
\]
From this identity, the strong continuity of the free group and the estimate \eqref{eq:energybd},
using the absolute continuity of the $L^2$ integral of the force, we obtain right continuity of the
state. Applying the same identity backwards with $S_\theta(-h)$ gives left continuity. We conclude
that $Z(t):=S_\theta(t)z+R_{\theta,t}p$ belongs to $C([0,T];\X)$.

\noindent Let us check that the displacement defined in this way satisfies the expected variational equation:
\begin{equation}\label{eq:weak}
 \frac{d}{dt}(u_t(t),v)_H+\form_\theta(u(t),v)=p(t)\cdot v(0),
 \qquad v\in V,\quad \text{in }\mathcal D'(0,T).
\end{equation}
Let $(e_n)_n$ be an $H$-orthonormal basis of real eigenfunctions of $\op_\theta$, with associated
eigenvalues $\omega_n^2>0$. Taking $h=(e_n,0)$ and $h=(0,e_n)$ in \eqref{eq:reach}, we obtain, with
$q_n=(u,e_n)_H$,
\[
 q_n(t)=q_n(0)\cos\omega_nt+
 q_n'(0)\frac{\sin\omega_nt}{\omega_n}
 +\int_0^t\frac{\sin(\omega_n(t-r))}{\omega_n}\,p(r)\cdot e_n(0)\,dr .
\]
The corresponding velocity coordinate is $q'_n$. By completeness of the eigenfunctions and
continuity of the state, the second component is indeed $u_t$. The equations for $q_n$ yield
\eqref{eq:weak} first for finite linear combinations of eigenfunctions and then for all $v\in V$
by density in the form norm. Conversely, \eqref{eq:weak}, together with the initial conditions,
determines all modal coordinates and hence the state uniquely. The same argument applies to weakly
continuous, energy-bounded states satisfying these conditions.

\noindent We shall also make use of the solution for concentrated interior sources. For a load of the form
$\sum_j g_j(t)b\,\delta_{\xi_j}$ on the right-hand side of the equation, we replace the integrand in
\eqref{eq:reach} by $\sum_j g_j(t)b\cdot\psi_t^h(\xi_j,t)$. By Lemma~\ref{lem:trace}, this defines a
bounded operator $P_{\theta,T}:L^2(0,T)^M\to\X$ giving the terminal state from rest with homogeneous
boundary conditions. Arguing as above by means of the modal equations, we obtain the weak equation
with right-hand side $\sum_j g_j(t)s_v(\xi_j)$.

\begin{lemma}\label{lem:compact}
For each fixed configuration $\theta$ and every $T>0$, $R_{\theta,T}-R_{\refn,T}:\U_T\to\X$ is compact.
\end{lemma}
\begin{proof}
Let $v$ and $u$ be the displacements of the uncoupled and coupled systems, respectively, associated
with the same force $p$ and zero initial data. Subtracting their variational equations and keeping
the coupled form on the left-hand side, we obtain
\[
 (R_{\theta,T}-R_{\refn,T})p=P_{\theta,T}g,\qquad g_j(t)=-k_j s_v(\xi_j,t).
\]
Since $P_{\theta,T}$ is bounded, it suffices to show that the map $p\mapsto g$ is compact. Consider a
bounded set of forces in $\U_T$. By \eqref{eq:energybd}, the corresponding $\|v(t)\|_V$ and
$\|v_t(t)\|_H$ are bounded uniformly in $t$. We use the following elementary estimate: for every
scalar function $w\in H^1(0,L)$ such that $w(L)=0$,
\begin{equation}\label{eq:interpol}
 |w(x)|^2\le 2\|w\|_{L^2}\|w'\|_{L^2}
\end{equation}
which follows by integrating $(w^2)'$ between $x$ and $L$. Applying \eqref{eq:interpol} to
$v_i(t)-v_i(r)$, whose $L^2$ norm is bounded by $C|t-r|$ and whose derivative is bounded in $L^2$,
we obtain equicontinuity in time in the uniform norm. Together with the spatial $1/2$-H\"older
estimate, this gives uniform boundedness and equicontinuity of the displacements on
$[0,L]\times[0,T]$. By the Arzel\`a--Ascoli theorem, from every such sequence we can extract a uniformly
convergent subsequence, and hence a subsequence for which the $g_j$ converge strongly in $L^2(0,T)$.
The conclusion follows from the boundedness of $P_{\theta,T}$.
\end{proof}

\section{Exact controllability with two boundary forces}\label{sec:two}
\subsection{Observability at the critical time}
\noindent \rev{Let $T>0$ be the observation time and set $c=\min_i c_i$. To recover both travelling
components of member $i$ at $(x,\tau)$ from the boundary at $x=0$, the characteristic identity
uses the boundary times $\tau-x/c_i$ and $\tau+x/c_i$. For these times to belong to $[0,T]$
for both members and all $0<x<L$, the middle slice must satisfy
\[
  \frac{L}{c}\leq \tau\leq T-\frac{L}{c}.
\]
Such a slice exists precisely when $T\geq 2L/c=T_*$. At equality it is unique,
with $\tau=L/c=T_*/2$. Set $Q_*=[0,L]\times[0,T_*]$ and
$I(x)=(x/c,T_*-x/c)$. Then $\tau\in I(x)$ for every $0\leq x<L$.
The open interval collapses at $x=L$, where the displacement is recovered by spatial continuity.}
We use the characteristic variables of the
one-dimensional wave equation (see~\cite{DagerZuazua}), taking into account their jumps at the
connectors.

\begin{lemma}\label{lem:side}
Let $(\psi,\psi_t)=S_\theta(\cdot)y$ and $g=O_{\theta,T_*}y$. For every finite configuration satisfying
$k_j\geq0$ and $\sum_j k_j\leq K$, we have
\begin{equation}\label{eq:weakobs}
 \|y\|_\theta^2\le C\left(\|g\|_{\U_{T_*}}^2+\|s_\psi\|_{C(Q_*)}^2\right),
\end{equation}
where $C$ depends only on $L,K$ and the material coefficients.
\end{lemma}
\begin{proof}
We first consider initial data given by finite linear combinations of eigenfunctions. Introduce the
characteristic variables $r_i^\sigma=\psi_{i,t}+\sigma c_i\psi_{i,x}$, $\sigma\in\{-1,1\}$. Between two
successive connectors, we have $(\partial_x-\sigma c_i^{-1}\partial_t)r_i^\sigma=0$. Since the velocity
is continuous across $\xi_j$, Proposition~\ref{prop:domain} yields the jump relation
$[r_i^\sigma]_{\xi_j}=\sigma(c_i/a_i)k_j b_i s_\psi(\xi_j,\cdot)$. The homogeneous Neumann condition
implies that $r_i^\sigma(0,\cdot)=g_i$. Integrating along the corresponding characteristic from
$x=0$ to $(x,t)$, we deduce, for $t\in I(x)$ and $x\notin\{\xi_1,\ldots,\xi_M\}$,
\begin{equation}\label{eq:chartransport}
 \begin{aligned}
 r_i^\sigma(x,t)
 &=g_i(t+\sigma x/c_i)\\
 &\quad+\sigma\frac{c_i b_i}{a_i}\sum_{\xi_j<x}k_j
 s_\psi\!\left(\xi_j,t+\sigma(x-\xi_j)/c_i\right).
 \end{aligned}
\end{equation}
Indeed, since $c_i\geq c$, we have $t+\sigma(x-z)/c_i\in I(z)$ whenever $0\leq z\leq x$ and
$t\in I(x)$. Thus the terms on the right-hand side of \eqref{eq:chartransport} involve only observed
boundary times and slip values at the connectors within the same characteristic region.

\noindent \rev{At $t=\tau$, the boundary arguments $\tau\pm x/c_i$ lie in $[0,T_*]$ for the entire
range $0<x<L$. The endpoint $x=L$ has zero measure in the energy integral. For $T<T_*$ the
same construction at $t=T/2$ would cover only $0<x<cT/2$ in the slower member. At $T=T_*$,
a change of variables gives}
\[
 \int_0^L |g_i(\tau+\sigma x/c_i)|^2\,dx
 \le c_i\int_0^{T_*}|g_i(r)|^2\,dr .
\]
The connector sum in \eqref{eq:chartransport} is bounded in absolute value by
$(c_i/a_i)K\|s_\psi\|_{C(Q_*)}$. By energy conservation and $a_i=\rho_i c_i^2$, we have
\[
 \|y\|_\theta^2
 =\frac12\sum_{i,\sigma}\rho_i\int_0^L|r_i^\sigma(x,\tau)|^2\,dx
   +\sum_j k_j|s_\psi(\xi_j,\tau)|^2
\]
Combining these estimates gives \eqref{eq:weakobs} for the regular initial data.

\noindent To pass to general initial data in $\X$, we approximate by finite linear combinations of
eigenfunctions. The corresponding states converge in $C([0,T_*];\X)$ and their displacements
converge uniformly on $Q_*$. By Lemma~\ref{lem:trace}, the observations converge in $\U_{T_*}$.
The preceding change of variables also ensures convergence of the translated observations
$g_i(\tau+\sigma x/c_i)$ in $L^2(0,L)$. We can therefore pass to the limit in both the identity on
the middle time slice and the estimate. Notice also that \eqref{eq:chartransport} is valid almost
everywhere in the open characteristic region: by integrating first in $t\in I(x)$ and then in
$x$, the squared norm of the translated observation is bounded by $L\|g_i\|_{L^2(0,T_*)}^2$.
\end{proof}

\begin{lemma}\label{lem:unique}
For every finite configuration with $k_j\geq0$, $O_{\theta,T_*}y=0$ implies $y=0$.
\end{lemma}
\begin{proof}
Assume that $g=O_{\theta,T_*}y=0$. By Lemma~\ref{lem:trace}, the distributional time derivative
of $\psi(0,\cdot)$ vanishes. Hence $\psi(0,t)=d$ for some constant vector $d$. We construct a
continuous, piecewise affine function $q$ by prescribing its value and derivative at $0$ and
continuing it successively across the connectors by means of the jump conditions:
\[
 q(0)=d,\qquad Dq'(0)=0,\qquad [Dq']_{\xi_j}=k_j b\,s_q(\xi_j).
\]
At this stage we do not impose the condition at $L$ on $q$; its value there is determined by the
above construction. Set $w(x,t)=\psi(x,t)-q(x)$. Since $q$ is independent of time and satisfies
both the equations between connectors and the jump conditions, \eqref{eq:chartransport} holds
for $w$ with vanishing boundary term and $s_w$ in place of $s_\psi$. This can be seen directly by
subtracting the characteristic identity for $q$, obtained by summing its slope jumps.

\noindent On the part of the characteristic region lying before the first connector, the sum on the
right-hand side is empty. Both characteristic variables of $w$ vanish and hence $w_t=w_x=0$ there.
As $w(0,t)=0$, we deduce that $w$ vanishes in this region. Suppose now that $w$ vanishes up to
the $j$-th connector. By continuity, $w(\xi_j,t)=0$ for $t\in I(\xi_j)$. Consequently, all the slip
terms in \eqref{eq:chartransport} vanish on the next strip. The characteristic variables again
vanish, and the constant value of $w$ on this strip is zero by continuity across the connector.
By induction, $w$ vanishes in the whole open characteristic region. \rev{At the critical time
this gives $\psi(x,\tau)=q(x)$ for every $x<L$. Taking $x\uparrow L$ and using spatial
continuity extends the identity to the fixed end.} From the characteristic identity for $w$ on the middle
time slice we also have $\psi_t(\cdot,\tau)=0$ in $H$.

\noindent \rev{Thus $q(L)=\psi(L,\tau)=0$. This step puts $q$ in $V$ and permits testing the
stationary equation with $q$.} By integrating by parts on the
subintervals and using the Neumann and transmission conditions for $q$, we have
$\form_\theta(q,v)=0$ for every $v\in V$. Choosing $v=q$, we obtain
\[
 0=\int_0^L Dq'\cdot q'\,dx+\sum_j k_j|s_q(\xi_j)|^2.
\]
Since all the stiffnesses are nonnegative and $q(L)=0$, we deduce that $q=0$. Thus the state
vanishes at time $\tau$ and, by reversibility of the free evolution, $y=0$.
\end{proof}

\noindent Let us first explain how the above uniqueness result, together with the compactness result of
Lemma~\ref{lem:compact}, yields exact controllability at $T_*$ for each fixed configuration. In fact,
for the uncoupled system \eqref{eq:chartransport} gives
\[
 \|y\|_{\refn}^2
 =\frac12\sum_i\rho_i c_i
       \int_{\tau-L/c_i}^{\tau+L/c_i}|g_i(r)|^2\,dr
 \le \frac12\max_i(\rho_i c_i)\,\|g\|_{\U_{T_*}}^2 .
\]
By \eqref{eq:reach} and the closed range theorem, the operator $R_{\refn,T_*}$ is surjective and
admits a bounded right inverse $B$. Consequently
$R_{\theta,T_*}B=I+(R_{\theta,T_*}-R_{\refn,T_*})B$ is the identity plus a compact operator.
By the Fredholm alternative, its range is closed and of finite codimension. Since the range of
$R_{\theta,T_*}$ contains this subspace, it is closed as well. On the other hand, it is dense by
Lemma~\ref{lem:unique} and \eqref{eq:reach} with $h=S_\theta(T_*)y$. This proves surjectivity.
For the optimization problem, however, an observability estimate for each fixed configuration is
not sufficient. We need an estimate with a constant which remains uniform as the configuration
varies, including when some of the stiffnesses vanish. This is the content of the next proposition.

\begin{proposition}\label{prop:unif}
Let $\mathcal K$ be a compact set of configurations with $M$ connectors of nonnegative stiffnesses
and bounded total stiffness. Assume that the distances between successive connectors, as well as
those from the outermost connectors to the boundary, are bounded from below by a positive constant.
Then there exists $\beta_*>0$, independent of $\theta\in\mathcal K$ and $T\geq T_*$, such that for every
$y\in\X$,
\begin{equation}\label{eq:unifobs}
 \int_0^T|\psi_t(0,t)|^2\,dt\ge\beta_*\|y\|_{\refn}^2,
 \qquad (\psi,\psi_t)=S_\theta(\cdot)y .
\end{equation}
\end{proposition}
\begin{proof}
It is enough to prove the estimate for $T=T_*$, since the observation integral is nondecreasing
with respect to $T$. We argue by contradiction. If the estimate does not hold, by compactness of
$\mathcal K$ there exist sequences $\theta_n\to\theta\in\mathcal K$ and $y_n\in\X$ such that
$\|y_n\|_{\refn}=1$ and $O_{\theta_n,T_*}y_n\to0$ in $\U_{T_*}$. Extracting a subsequence, we may assume
that $y_n\rightharpoonup y$ in $\X$. Set $\Psi_n(t)=S_{\theta_n}(t)y_n$. For each fixed $t$ and $h\in\X$,
we have
\[
 (\Psi_n(t),h)_{\theta_n}
 =(y_n,S_{\theta_n}(-t)h)_{\theta_n}
 \longrightarrow (y,S_\theta(-t)h)_\theta .
\]
In deriving this convergence we used \eqref{eq:moving}, which gives convergence of the energy
scalar products in form norm, and the strong convergence of the groups on each fixed datum
established in Section~\ref{sec:forms}. Since the states are uniformly bounded by \eqref{eq:uniform},
we can replace the scalar product on the left by $(\cdot,h)_\theta$ when taking the limit. Hence
$\Psi_n(t)\rightharpoonup S_\theta(t)y$ in $\X$ for every $t$.

\noindent We next establish strong convergence of the displacements. The uniform energy bound gives
spatial equicontinuity. In addition, $\|\psi_n(t)-\psi_n(r)\|_H\leq C|t-r|$, which, together with
\eqref{eq:interpol}, yields equicontinuity in time in the uniform spatial norm. By the
Arzel\`a--Ascoli theorem, up to a subsequence the displacements converge uniformly on $Q_*$.
By the weak convergence at every time and the continuity of point evaluation on $V$, the limit
is the displacement $\psi$ associated with $S_\theta(\cdot)y$. Since the observations converge to
zero in $L^2$, the limit trace $\psi(0,\cdot)$ has zero distributional time derivative. It follows
from Lemma~\ref{lem:unique} that $y=0$. In particular, $\|s_{\psi_n}\|_{C(Q_*)}\to0$. This is precisely
the convergence needed to estimate the slip terms in \eqref{eq:weakobs} when the connector positions
vary. Applying that estimate to $(\theta_n,y_n)$, we deduce that $\|y_n\|_{\theta_n}\to0$, in
contradiction with $1=\|y_n\|_{\refn}\leq\|y_n\|_{\theta_n}$.
\end{proof}

\subsection{The control of minimum norm}
\noindent We derive here the characterization of the minimum-norm control, which will be used in the design
problem. Given the initial state $z$ and the target state $z^T$, set $w=S_\theta(-T)z^T-z$. Taking
$h=S_\theta(T)y$ in \eqref{eq:reach}, we obtain, for every $T>0$,
\begin{equation}\label{eq:dual}
 Z(T)=z^T\quad\Longleftrightarrow\quad
 (p,O_{\theta,T}y)_{\U_T}=(w,y)_\theta\quad(y\in\X).
\end{equation}
Assume that $T\geq T_*$. By the Hilbert Uniqueness Method (see~\cite{Lions88,TucsnakWeiss}), we are led
to the following minimization problem and the associated control:
\begin{equation}\label{eq:HUM}
 J_w(y)=\tfrac12\|O_{\theta,T}y\|_{\U_T}^2-(w,y)_\theta,\qquad
       p_{\min}=O_{\theta,T}\widehat y .
\end{equation}
Here $\widehat y$ denotes the minimizer of $J_w$ over $\X$. Indeed, by \eqref{eq:unifobs}, the continuous
bilinear form $(O_{\theta,T}\cdot,O_{\theta,T}\cdot)_{\U_T}$ is coercive on $\X$. By the Lax--Milgram
theorem, there exists a unique $\widehat y$ solving
$(O_{\theta,T}\widehat y,O_{\theta,T}y)_{\U_T}=(w,y)_\theta$ for every $y\in\X$, which is also the unique
minimizer of $J_w$. It follows from \eqref{eq:dual} that $p_{\min}$ steers the system to $z^T$. The set of
all controls reaching this target is $p_{\min}+(\operatorname{Ran}O_{\theta,T})^\perp$. As $p_{\min}$ belongs
to $\operatorname{Ran}O_{\theta,T}$, it is the unique control of minimum norm. In particular, for the
problem of bringing the system to rest, we have
\begin{equation}\label{eq:costdual}
 m_T(\theta,z)=\sup_{y\ne0}\frac{|(z,y)_\theta|}{\|O_{\theta,T}y\|_{\U_T}}.
\end{equation}
To see this, take $w=-z$ in \eqref{eq:dual}. By the Cauchy--Schwarz inequality, the quotient is
bounded by $\|p_{\min}\|_{\U_T}$. If $z\neq0$, equality is attained at $y=\widehat y$, whereas for $z=0$
both sides vanish.

\subsection{Necessity of the control time}
\noindent To prove the necessity of the condition $T\geq T_*$, we exploit the fact that the endpoint map of
the coupled system differs from that of the uncoupled one by a compact operator. Assume that
$0<T<T_*$ and consider a member whose wave speed is $c=\min_i c_i$. Set $P=2L/c$. We first construct
an infinite-dimensional family of free solutions of the uncoupled system which are not observed
at the accessible end during $(0,T)$. Given any $q\in C_c^\infty(T,P)$, extend it to $\R$ by imposing
$q(t+P)=-q(t)$ and set
\[
 F(t)=-\tfrac12\int_0^P q(r)\,dr+\int_0^t q(r)\,dr .
\]
One checks that $F(t+P)=-F(t)$ by differentiation and evaluation at $t=0$. Consider the wave
$v(x,t)=F(t+x/c)+F(t-x/c)$ in the slow member and set the other component to zero. We have
$v_x(0,t)=0$, $v(L,t)=0$ and $v_t(0,t)=2q(t)=0$ for $t\in(0,T)$. The map from $q$ to the state of this
solution at time $T$ is injective: if the terminal state is zero, then the free solution vanishes
for all times and hence so does $q$. By \eqref{eq:reach}, these terminal states constitute an
infinite-dimensional subspace of $(\operatorname{Ran}R_{\refn,T})^{\perp_{\refn}}$.

\noindent Let $(h_n)_n$ be an orthonormal sequence in this orthogonal complement, with respect to the
reference scalar product. In the rest of the argument all adjoints are taken with respect to this
same scalar product on $\X$. Since $h_n\rightharpoonup0$ in $\X$ and the adjoint of the operator in
Lemma~\ref{lem:compact} is compact, $(R_{\theta,T}-R_{\refn,T})^*h_n\to0$ in $\U_T$. As
$R_{\refn,T}^*h_n=0$, we deduce that $R_{\theta,T}^*h_n\to0$ in $\U_T$. In view of $\|h_n\|_{\refn}=1$,
this contradicts the lower bound on the adjoint which would follow from surjectivity of
$R_{\theta,T}$. Hence the system cannot be exactly controllable for $T<T_*$, and the proof of the
first assertion of Theorem~\ref{thm:control} is complete.

\section{The case of one controlled member}\label{sec:one}
\noindent In contrast to systems with controls applied at the junctions (see, for instance, Ho~\cite{Ho}),
the connectors in our model are passive. We assume in this section that the force is applied
only to the first member at $x=0$, i.e., $p=(p_1,0)$. Denote by $\pi_2$ the projection of the state
onto the displacement and velocity of the second member. As the members are independent in the
uncoupled system, $\pi_2R_{\refn,T}(p_1,0)=0$. Hence, by Lemma~\ref{lem:compact}, the map
$p_1\mapsto\pi_2 R_{\theta,T}(p_1,0)$ is compact from $L^2(0,T)$ into
$\{v\in H^1(0,L):v(L)=0\}\times L^2(0,L)$. Since this space is infinite-dimensional, the map
cannot be surjective. Thus the system is not exactly controllable in $\X$ by one boundary force.

\noindent \rev{If all stiffnesses vanish,
the second member stays at rest. Otherwise let $\xi_a$ be the first connector of positive
stiffness, $k=k_a$ and $t_a=\xi_a/c_1$. Take the solution of the uncoupled system driven from
rest by $(p_1,0)$. For $t<t_a$, finite propagation at speed $c_1$ gives zero displacement on
member~1 at every active connector, while member~2 is at rest. Thus all connector forces vanish
on this interval, and the uncoupled solution also solves \eqref{eq:weak} for $\theta$. By
uniqueness it agrees there with the coupled solution. Continuity in the energy space gives
$K_{\theta,T}=0$ also at $T=t_a$.}

\noindent \rev{Suppose $T>t_a$ and put $c_{\max}=\max_i c_i$. Choose $0<h<T-t_a$ so small that
$c_{\max}h$ is less than the distance from $\xi_a$ to either endpoint and to every other
active connector. Start the force at $t_0=T-h-t_a$. Its first possible arrival at $\xi_a$ is
then $t_0+t_a=T-h$. During the remaining interval of length $h$ only this connector is excited,
and the waves emitted there cannot reach an endpoint or another active connector. Write
$Z_i=\rho_i c_i=a_i/c_i$ for the characteristic impedance and
$s(t)=u_1(\xi_a,t)-u_2(\xi_a,t)$ for the slip. The incident right-going velocity on member~1
at $\xi_a$ is $p_1(t-t_a)/Z_1$, by the boundary condition $-a_1u_{1,x}(0,t)=p_1(t)$.}

\noindent \rev{Subtract the freely propagating incident wave on member~1 and denote the scattered
displacement on each member by $w_i$ (thus $w_2=u_2$). On the time interval above, $w_i$
consists of a left-going wave to the left of $\xi_a$ and a right-going wave to its right.
Continuity gives a common scattered velocity trace $v_i^{\rm out}=\partial_t w_i(\xi_a,t)$.
Hence $a_i w_{i,x}(\xi_a^-,t)=Z_i v_i^{\rm out}$ and
$a_i w_{i,x}(\xi_a^+,t)=-Z_i v_i^{\rm out}$. The incident field has continuous axial force,
so the point force $F_i$ gives $[a_i w_{i,x}]=-2Z_i v_i^{\rm out}=-F_i$, and
$v_i^{\rm out}=F_i/(2Z_i)$. Here $F_1=-ks$ and $F_2=ks$. Subtracting the two velocities at
the connector gives
\[
 s'(t)+\frac{k}{2}\bigl(Z_1^{-1}+Z_2^{-1}\bigr)s(t)
    = Z_1^{-1}p_1(t-t_a),\qquad s(T-h)=0 .
\]
Consequently, each $f\in C_c^\infty(T-h,T)$ can be realized as the slip on $(T-h,T)$.
Choose $p_1(t-t_a)=Z_1[f'(t)+\tfrac{k}{2}(Z_1^{-1}+Z_2^{-1})f(t)]$ there and extend the
force by zero outside the corresponding input interval. Transport of the outgoing wave to the
right on member~2 gives
\[
 u_{2,t}(x,T)=\frac{k}{2Z_2}\,
       f\!\left(T-\frac{x-\xi_a}{c_2}\right),
 \qquad \xi_a<x<\xi_a+c_2h .
\]
Restriction of the terminal velocity to this segment is injective in $f$. Hence the projected
endpoint map has infinite rank. If its range were closed, the restriction to the orthogonal
complement of the kernel would have a bounded inverse. Composing with that inverse would make
the identity on the infinite-dimensional range compact. This proves the second assertion of
Theorem~\ref{thm:control}.}

\begin{remark*}
\rev{Write $\X_2=\{v\in H^1(0,L):v(L)=0\}\times L^2(0,L)$ with its scalar energy inner
product, whose squared norm is $\int_0^L(a_2|v'|^2+\rho_2|w|^2)\,dx$ for $(v,w)\in\X_2$.
For $K=K_{\theta,T}$, let $(\sigma_n,e_n,f_n)$ be its singular system for the positive singular
values, with $Ke_n=\sigma_n f_n$. By the Picard criterion the reachable targets of the second
member satisfy the equivalence
\[
 z_2\in\operatorname{Ran}K
 \quad\Longleftrightarrow\quad
 z_2\perp\ker K^*
 \quad\text{and}\quad
 \sum_n\frac{|(z_2,f_n)_{\X_2}|^2}{\sigma_n^2}<\infty .
\]
Indeed, the $f_n$ form an orthonormal basis of the closure of the range, and the coefficients of a minimum-norm
force in the $e_n$ are $(z_2,f_n)_{\X_2}/\sigma_n$. Orthogonality alone characterizes approximate
reachability in $\X_2$. These conditions concern the projection onto member~2.}
\end{remark*}

\noindent \rev{For the full one-force endpoint map $W_Tp_1=R_{\theta,T}(p_1,0)$, the closure in
$\X$ is $S_\theta(T)(\ker O_{\theta,T}^{(1)})^{\perp_\theta}$, where orthogonality is in the
coupled energy scalar product.} By \eqref{eq:reach}, the density
of the reachable set in $\X$ is equivalent to injectivity of the scalar observation
$O_{\theta,T}^{(1)}y=(O_{\theta,T}y)_1$. Two necessary conditions follow immediately: every eigenspace of
$\op_\theta$ must be one-dimensional and every nonzero eigenfunction $e$ must satisfy $e_1(0)\neq0$.
Indeed, for $\op_\theta e=\omega^2e$, the free solution $\psi(x,t)=e(x)\sin(\omega t)/\omega$ is not
observed if $e_1(0)=0$. Moreover, on an eigenspace of dimension at least two, the linear functional
$e\mapsto e_1(0)$ has a nontrivial kernel. Note that the dimension of each eigenspace is at most two:
the two values of the displacements at $0$, together with the homogeneous Neumann condition,
determine the solution of the system of ODEs successively across the connectors. We do not claim
that the above necessary conditions are sufficient for approximate controllability in $\X$.
However, they are sufficient for the corresponding problem in a finite spectral subspace, as
shown in Section~\ref{sec:modal}. We give below two examples illustrating the dependence of the
invisible modes on the parameters of the system.

\begin{example}\label{ex:node}
Consider a single connector at $\xi=L/3$. The function defined by
\[
 e_1=0,\qquad e_2(x)=\cos(3\pi x/(2L)),\qquad
 \lambda=c_2^2(3\pi/(2L))^2
\]
satisfies the homogeneous boundary conditions and $\op_\theta e=\lambda e$ for every $k\geq0$.
Indeed, the slip $s_e(\xi)$ vanishes and the axial forces are continuous across the connector.
This mode cannot be controlled by a force applied to the first member. More generally, the same
obstruction occurs if all connectors of positive stiffness are located at zeros of the same
function $\cos((2m+1)\pi x/(2L))$, $m\geq0$. The slip for this mode vanishes at each active connector,
independently of its stiffness. In particular, increasing the stiffnesses at these positions
cannot restore the controllability of this mode.
\end{example}

\begin{example}\label{ex:ratio}
The vanishing of the slip at the connectors is not necessary for a mode to be invisible. To see
this, consider one connector $(\xi,k)$ with $k>0$ and set $\kappa_i=\sqrt{\lambda}/c_i$. An invisible
eigenfunction, normalized by $e_2(0)=1$, satisfies $e_1=0$ and $e_2=\cos(\kappa_2x)$ on $(0,\xi)$.
Continuing it to $x>\xi$ by the transmission conditions, we obtain
\[
 \begin{aligned}
 e_1(x)&=-\frac{k\cos(\kappa_2\xi)}{a_1\kappa_1}\sin(\kappa_1(x-\xi)),\\
 e_2(x)&=\cos(\kappa_2 x)
       +\frac{k\cos(\kappa_2\xi)}{a_2\kappa_2}\sin(\kappa_2(x-\xi)).
 \end{aligned}
\]
Assume that $\cos(\kappa_2\xi)\neq0$. The condition $e_1(L)=0$ implies that
$\kappa_1=n\pi/(L-\xi)$, $n\geq1$. Substituting $\kappa_2=(c_1/c_2)\kappa_1$ into $e_2(L)=0$ gives an
equation for the ratio of the two wave speeds. For instance, take $L=1$, $\xi=1/2$, $k=a_2=1$ and $n=1$,
and set $s=2\pi c_1/c_2$. The equation is $\cos s+\sin s/(2s)=0$. On each interval
$((m+\frac12)\pi,(m+1)\pi)$, $m\geq0$, the function $\tan s+2s$ is strictly increasing from $-\infty$
to a positive value. Its unique zero satisfies $\cos(s/2)\neq0$. The first four corresponding ratios
$c_1/c_2$ are approximately $0.29230$, $0.76647$, $1.26004$ and $1.75720$. Thus, even when the connector
exerts a nonzero force, an eigenmode can have zero trace at the controlled end of the first member.
The absence of connectors at vibration nodes is therefore not, on its own, a sufficient design
condition for control by a single force.
\end{example}
\noindent The examples illustrate why the choice of the positions and stiffnesses matters in the one-force
problem, despite the impossibility of exact controllability in the energy space. To obtain control
results in weaker, spectrally weighted spaces, one would have to study quantitative estimates of
the observed traces and the spectral separation. In the present work, we restrict this part of the
analysis to a prescribed finite set of modes.

\section{Optimal configurations of connectors in the presence of damage}\label{sec:design}
\subsection{Dependence of the observability constant on the configuration}
\noindent We turn now to the optimization problem formulated in \eqref{eq:design}. Recall that we have fixed
one and the same set of initial states, $\{z:\|z\|_{\refn}=1\}$, for all the configurations. This is
important in comparing the costs for different configurations, since the energy scalar product in
\eqref{eq:costdual} depends on the stiffness measure. The relevant dual norm is given by
\[
 D_\theta(y)=\sup_{\|z\|_{\refn}=1}|(z,y)_\theta|
            =\|L_\theta y\|_{\refn},\qquad
 (z,L_\theta y)_{\refn}=(z,y)_\theta .
\]
Here $L_\theta$ denotes the Riesz operator on the reference space. By \eqref{eq:uniform}, we have
$I\leq L_\theta\leq(1+K\alpha^2L)I$ for configurations whose total stiffness does not exceed $K$. In
particular, $\|y\|_{\refn}\leq D_\theta(y)\leq(1+K\alpha^2L)\|y\|_{\refn}$. From \eqref{eq:moving} we also
have continuity of $L_\theta$ with respect to $\theta$ in the operator norm. Taking the supremum of
\eqref{eq:costdual} over $\|z\|_{\refn}=1$ and interchanging the two suprema, we obtain
\begin{equation}\label{eq:best}
 \mathcal C_T(\theta)=\beta_T(\theta)^{-1/2},\qquad
 \beta_T(\theta)=\inf_{y\ne0}
        \frac{\|O_{\theta,T}y\|_{\U_T}^2}{D_\theta(y)^2}.
\end{equation}
For $T\geq T_*$, Proposition~\ref{prop:unif}, combined with the upper bound on $D_\theta(y)$, yields
$\beta_T(\theta)\geq\beta_*/(1+K\alpha^2 L)^2$ for the configurations in a compact admissible family.
This estimate bounds the worst control cost independently of the configuration. To prove that the
minimum is attained, we establish first a continuity property of the observations.

\begin{lemma}\label{lem:strongtrace}
Let $\theta_n\to\theta$ in a compact admissible family with a common positive boundary clearance.
Then, for every $T>0$ and every fixed $y\in\X$, we have
$\|O_{\theta_n,T}y-O_{\theta,T}y\|_{\U_T}\to0$.
\end{lemma}
\begin{proof}
By the result of Section~\ref{sec:forms},
$\sup_{0\leq t\leq T}\|(S_{\theta_n}(t)-S_\theta(t))y\|_{\refn}\to0$. Let $w$ denote the difference of the
corresponding displacements. The interval $(0,d/2)$ contains no connector for any of the
configurations. Hence $w$ satisfies the two free scalar equations in this interval, with the
homogeneous Neumann condition at $0$. Applying \eqref{eq:mult} with $q(0)=1$ and $q(d/2)=0$, we bound
its right-hand side by $C(1+T)\sup_t\|(w,w_t)(t)\|_{\refn}^2$. It follows that the difference of the
observations tends to zero in $L^2(0,T)^2$. To justify the multiplier identity for energy solutions,
we approximate each wave by finite sums of its own eigenfunctions. Their difference satisfies the
same free equations on the common interval, the states converge uniformly in the energy norm and
the boundary observations converge by Lemma~\ref{lem:trace}. This allows us to pass to the limit.
Note that we have used strong convergence for each fixed datum, not convergence in the operator
norm of the observation operators.
\end{proof}

\begin{proof}
\rev{To prove Theorem~\ref{thm:design}, fix $T\geq T_*$. The set $\Gamma_r$ is a finite union of compact sets, obtained by
fixing all coordinates outside a subset of at most $r$ indices equal to one. The family
\[
 \mathcal K=\{\theta^\gamma:\theta\in\Tad,\ \gamma\in\Gamma_r\}
\]
is the continuous image of the compact product $\Tad\times\Gamma_r$.} In addition, the total stiffness of its elements is bounded by $K$, and their
positions satisfy the original clearance and separation conditions. Hence, by
Proposition~\ref{prop:unif} and \eqref{eq:best}, the cost $\mathcal C_T$ is bounded from above on
$\mathcal K$, including configurations with vanishing stiffnesses.

\noindent For each $y\in\X\setminus\{0\}$, the map
$\theta\mapsto D_\theta(y)/\|O_{\theta,T}y\|_{\U_T}$ is continuous on $\mathcal K$, by
Lemma~\ref{lem:strongtrace}, estimate \eqref{eq:moving} and the positivity of the denominator.
Consequently, its supremum over $y\neq0$, which is precisely $\mathcal C_T(\theta)$, is lower
semicontinuous. \rev{For each fixed $\gamma\in\Gamma_r$, the map
$\theta\mapsto\mathcal C_T(\theta^\gamma)$ is lower semicontinuous on $\Tad$.
Taking the supremum over $\Gamma_r$ gives a lower semicontinuous functional
$\mathcal J_{T,r}$, finite by the uniform bound above.} Let $(\theta_n)_n$ be a minimizing sequence. By compactness, up to a subsequence,
$\theta_n\to\theta_*\in\Tad$. We have
\[
 \rev{\mathcal J_{T,r}(\theta_*)\leq\liminf_{n\to\infty}\mathcal J_{T,r}(\theta_n)
   =\inf_{\theta\in\Tad}\mathcal J_{T,r}(\theta),}
\]
which proves the assertion. Note that neither the attainment of the supremum over the initial
states nor that of the \rev{supremum over $\Gamma_r$} is needed in this argument.
\end{proof}

\noindent Let us emphasize the order of the optimization in \eqref{eq:design}. The nominal configuration is
common to all damage scenarios. For each of these scenarios and each initial state, the force is
chosen by solving the corresponding control problem. Requiring one and the same force to work for
several unknown damaged configurations would be a simultaneous controllability problem, which is
different from the one addressed here. Also, lower semicontinuity suffices for the existence of
an optimal nominal configuration, but does not by itself guarantee the existence of \rev{a worst damage scenario}.

\noindent \rev{The intact state also fits the supremum in \eqref{eq:design} when the losses can be
arbitrarily small. Indeed, if $\gamma_0<1$, take $\gamma^{(n)}\in\Gamma_r$ with one
coordinate strictly below one and tending to one, and all others equal to one. Then
$\mathcal C_T(\theta)\leq\liminf_n\mathcal C_T(\theta^{\gamma^{(n)}})$ by lower
semicontinuity. Including the intact vector explicitly therefore gives the same supremum.}

\begin{example}\label{ex:freecost}
\rev{Take $r=M$ and $\gamma_0=0$. Complete failure of all nominal connectors is admissible
for every $\theta\in\Tad$, since $\gamma=(0,\ldots,0)$ belongs to $\Gamma_M$.} The resulting system consists
of two independent beams, and the optimal observability constant can be evaluated explicitly.
Set $P_i=2L/c_i$. For every $T>0$, the infimum in \eqref{eq:best} is given by
\begin{equation}\label{eq:freebeta}
 \beta_T(\refn)=\min_{i=1,2}\frac{2\lfloor T/P_i\rfloor}{\rho_i c_i}.
\end{equation}
Indeed, a scalar free solution has the representation $v(x,t)=F(t+x/c_i)+F(t-x/c_i)$, where
$F(t+P_i)=-F(t)$. Denote by $q=F'$. This representation includes all the energy solutions: given
initial data, we prescribe $q(\pm x/c_i)=(v_t(x,0)\pm c_i v_x(x,0))/2$ for $0<x<L$ and extend $q$
antiperiodically. The integration constant in $F$ is determined by the antiperiodicity of $F$.
The resulting initial displacement has the prescribed derivative and vanishes at $L$. We have
for the norm and the boundary observation of this component
\[
 \|(v,v_t)(t)\|_{\refn,i}^2
   =2\rho_i c_i\int_{t-L/c_i}^{t+L/c_i}|q(r)|^2\,dr,\qquad
 \int_0^T|v_t(0,t)|^2\,dt=4\int_0^T|q(r)|^2\,dr .
\]
The first identity is obtained by expanding the squares and making a change of variables. Since
$|q|^2$ is periodic of period $P_i$, the integral of $|q|^2$ over any full period is the same. Write
$T=mP_i+r$ with $0\leq r<P_i$. The infimum of the ratio of the two integrals is realized by choosing
$q$ nonzero and supported in $(r,P_i)$ on one period, then extending it antiperiodically, so that
the contribution of the remaining interval of length $r$ vanishes. This includes the cases $m=0$
and $r=0$. The two components being independent, we obtain \eqref{eq:freebeta}.

\noindent \rev{For $r=M$, $\gamma_0=0$ and $T\geq T_*$, every nominal configuration therefore satisfies
\[
 \mathcal J_{T,M}(\theta)\geq\mathcal C_T(\refn)
 =\left(\min_{i=1,2}\frac{2\lfloor T/P_i\rfloor}{\rho_i c_i}\right)^{-1/2}.
\]
At $T=T_*$ each member contributes at least one complete period of observation. The lower
bound depends only on the uncoupled members and the time horizon, for any number and placement
of nominal connectors. Whether a design attains the bound requires information about the
remaining damage scenarios.}
\end{example}

\subsection{Control in a prescribed frequency band by one boundary force}\label{sec:modal}
\noindent In view of the compactness obstruction established in Section~\ref{sec:one}, it is natural to
formulate a more restricted control objective when we can act on only one of the members.
We consider the problem of eliminating the vibrations associated with a finite set of modes.

\noindent For the moment, fix the configuration $\theta$. Let $\{e_n:n\in I\}$ be a finite $H$-orthonormal
collection of eigenfunctions, $\op_\theta e_n=\omega_n^2e_n$, containing bases of the complete
eigenspaces of the selected frequencies. Denote by $Y_\theta$ the subspace of $\X$ spanned by
$h_n^d=(e_n/\omega_n,0)$ and $h_n^v=(0,e_n)$. These vectors form an orthonormal basis of $Y_\theta$
for the energy scalar product $(\cdot,\cdot)_\theta$. We require the orthogonal projection (in this
scalar product) of the terminal state $Z(T)$ onto $Y_\theta$ to vanish; no condition is imposed on the
other modal coordinates. As $Y_\theta$ is invariant under $S_\theta(t)$, the restriction of
\eqref{eq:dual} to $Y_\theta$ gives a necessary and sufficient condition for this problem. The scalar
observations of the basis vectors are given by
\begin{equation}\label{eq:modalobs}
 g_{n,d}(t)=-e_{n,1}(0)\sin(\omega_n t),\qquad
 g_{n,v}(t)= e_{n,1}(0)\cos(\omega_n t).
\end{equation}
Index the $2|I|$ basis vectors by $\ell$, and introduce the observation Gram matrix $G$ and the
matrix $H^{\refn}$ associated with the dual reference norm:
\[
 G_{\ell m}=\int_0^T g_\ell(t)g_m(t)\,dt,\qquad
 H_{\ell m}^{\refn}=(L_\theta h_\ell,L_\theta h_m)_{\refn}.
\]
\begin{proposition}\label{prop:band}
Let $T>0$. The system is controllable by one boundary force in the above sense, i.e., its terminal
projection on $Y_\theta$ can be made zero for every initial state, if and only if $G$ is positive
definite. Equivalently, each selected frequency is simple and $e_{n,1}(0)\neq0$ for every selected
mode. If these conditions hold, the worst control cost over $\|z\|_{\refn}=1$ is given by
\[
 \left(\min_{v\in\R^{2|I|}\setminus\{0\}}
       \frac{v^{\mathsf t}Gv}{v^{\mathsf t}H^{\refn}v}\right)^{-1/2}.
\]
If $G$ is singular, the worst control cost is infinite (with the convention that the cost for an
unattainable terminal projection is infinite).
\end{proposition}
\begin{proof}
For a given initial state $z$, the control constraints read
$\int_0^T p_1(t)g_\ell(t)\,dt=r_\ell$, where $r_\ell=-(z,h_\ell)_\theta$. Every vector $r$ can be obtained
by a suitable choice of $z\in Y_\theta$. If $G$ is positive definite, the control of minimum norm is
$\sum_\ell(G^{-1}r)_\ell g_\ell$, as is seen by taking its scalar products with $g_m$ and using its
orthogonality to the common annihilator of the $g_m$. If $G$ is singular, a nonzero vector in its
kernel gives a linear dependence between the constraints, which cannot be satisfied for some $r$.
This proves the characterization in terms of $G$.

\noindent On an eigenspace of dimension greater than one there exists a nonzero linear combination of the
eigenfunctions with zero first boundary trace. On the other hand, for distinct positive
frequencies, the functions $e^{\pm i\omega_n t}$ are linearly independent on every interval of
positive length. Indeed, differentiating a vanishing linear combination at an interior point
leads to a Vandermonde matrix associated with the distinct numbers $\pm i\omega_n$. In view of
\eqref{eq:modalobs}, this gives the equivalent spectral conditions. Finally, the same duality argument as in \eqref{eq:costdual}, now on $Y_\theta$ and with the scalar
observation, applies for every $T>0$ when $G>0$. For $y=\sum_\ell v_\ell h_\ell$, the squared observation
norm is $v^{\mathsf t}Gv$ and the squared dual reference norm is $v^{\mathsf t}H^{\refn}v$. Taking the
supremum first over the reference unit sphere and then over $y\in Y_\theta\setminus\{0\}$ gives the
stated generalized Rayleigh quotient.
\end{proof}

\noindent We can now formulate a design problem for the control of the modes in a prescribed frequency
band. Let $[\omega_-,\omega_+]\subset(0,\infty)$ be such that neither endpoint belongs to the spectrum
of $\sqrt{\op_\theta}$ for any configuration in the damaged family. We also assume that the number
of eigenfrequencies in this interval, counted with multiplicity, is fixed and positive. Denote by
$\beta_{\mathrm{band}}(\theta)$ the smallest generalized eigenvalue in Proposition~\ref{prop:band},
including the value zero if the matrix $G$ is singular. We claim that $\beta_{\mathrm{band}}$ depends
continuously on the configuration, even in the presence of eigenvalue crossings within the band.

\noindent To see this, fix a configuration and choose a contour enclosing precisely the eigenvalues of
$\op_\theta$ lying between the squared endpoints of the band. In a neighbourhood of this
configuration the same contour can be used. From the resolvent estimate associated with
\eqref{eq:moving}, we have convergence of $(\op_\theta+1)^{-1}$ in the operator norm from $H$ to $V$.
For $\zeta$ on the contour, we use the identity
\[
 (\op_\theta-\zeta)^{-1}
   =(\op_\theta+1)^{-1}\big[I-(1+\zeta)(\op_\theta+1)^{-1}\big]^{-1}.
\]
By continuity of inversion in the operator norm on $H$, uniformly on the contour near the fixed
configuration, the resolvents on the contour are continuous in $H\to V$. Integrating along the
contour, we deduce the continuity of the associated spectral projection in $H\to V$. Applying the
projection to a basis of the selected spectral subspace at the fixed configuration gives a
continuous local basis in $V$, and hence also in $H$. Orthonormalizing the corresponding phase
vectors for the energy scalar product gives a continuous local frame of $Y_\theta$. In this frame,
Lemma~\ref{lem:strongtrace} and the uniform direct estimate of Lemma~\ref{lem:trace} yield continuity
of the observed vectors in $L^2(0,T)$. Indeed, the variation due to the varying frame vectors is
bounded by the uniform direct estimate, and that due to the varying observation operator acting
on a fixed vector tends to zero by strong continuity. Thus the matrices $G$ and $H^{\refn}$ are
continuous in the local frame, with $H^{\refn}$ positive definite. The minimum of the Rayleigh
quotient on the Euclidean unit sphere is continuous. As this minimum is independent of the choice
of the frame, the claim follows.

\noindent By continuity on the compact set of design and damage parameters, the maximum in
\begin{equation}\label{eq:robust-band}
 \rev{\max_{\theta\in\Tad}\ \min_{\gamma\in\Gamma_r}
          \beta_{\mathrm{band}}(\theta^\gamma)}
\end{equation}
is attained. A strictly positive value of this maximum is equivalent to the existence of a nominal
configuration for which the modes in the prescribed band can be eliminated by one force, with a
uniform bound on the cost for all the damage scenarios. If, on the contrary, for each nominal
configuration there is a scenario in which at least one of the selected modes is invisible, the
maximum is zero. Example~\ref{ex:node} exhibits such a loss of observability for a particular
configuration and a band containing the indicated frequency.

\noindent \rev{If $r=M$ and $\gamma_0=0$, the robust value in \eqref{eq:robust-band} is zero
whenever the band contains a frequency $\omega_{2,n}=c_2(2n+1)\pi/(2L)$ of the uncoupled
second member, $n\geq0$. Indeed, $\gamma=(0,\ldots,0)$ is admissible for every nominal
design. In this scenario the selected eigenspace contains the mode
$e(x)=(0,\cos((2n+1)\pi x/(2L)))$. Its first boundary trace vanishes, so
Proposition~\ref{prop:band} gives $\beta_{\mathrm{band}}(\theta^\gamma)=0$.
The inner minimum is therefore zero for every $\theta\in\Tad$.}

\noindent The interest of \eqref{eq:robust-band} from the point of view of the design is that it distinguishes
between the loss of controllability of a mode and the deterioration of its control cost. Avoiding
vanishing boundary traces and repeated frequencies is necessary for controllability, but a small
positive generalized eigenvalue still entails a large worst-case force norm. It is the entire
Gram matrix, rather than only the nonvanishing of the individual traces, which quantifies the
quality of a configuration at the prescribed time. The maximization problem provides a criterion
for selecting the positions and the stiffnesses adapted to the band and to the damage scenarios.
Its positivity has to be investigated for the particular admissible class and the chosen band.
In particular, the projected control leaves the modes outside the band unconstrained and does not
imply that the full system is at rest.

\subsection{Further remarks on the design problem}\label{sec:scope}
\noindent Let us return to the mechanical motivation of this work. The results above concern the design
of a connection in view of the active control of longitudinal vibrations. In the case of two
boundary forces the minimum time is dictated by the wave speeds and not by the connectors. What
the optimization of the latter provides is a bound on the amount of boundary force needed to
eliminate an arbitrary initial vibration at the prescribed time, in the presence of damage. In
the case of one boundary force, the optimization has to be interpreted in a more restricted
sense, for instance in terms of the selected vibration modes. The two problems should not be
confused.

\noindent There is also a distinction between the active control criterion and a passive design criterion.
In the latter, given a bounded class of external forces $\mathcal P\subset\U_T$ applied to the
system initially at rest, one could, for instance, consider the worst slip
\[
 \rev{\sup_{\gamma\in\Gamma_r}\ \sup_{p\in\mathcal P}
     \|s_{u[\theta^\gamma,p]}\|_{C([0,L]\times[0,T])}.}
\]
Here \rev{$u[\theta^\gamma,p]$} denotes the displacement of the damaged system associated with the
load $p$ and zero initial data. This quantity is finite by \eqref{eq:energybd} and \eqref{eq:holder}.
It could be used, together with the \rev{connector force $|\gamma_j k_j s_u(\xi_j,t)|$ at each connector}, to assess the relative displacement of the two members and the demand on the
connection under a specified class of loads. Minimizing this criterion is a different problem
from minimizing the control cost in \eqref{eq:design}. In particular, our results do not imply that
an optimal configuration for one of these criteria is optimal for the other. Moreover, as the
connectors considered here are purely elastic, the unforced system conserves energy. To study
the decay of the vibrations one would have to introduce dissipation, for instance by means of
damping or of a feedback law.

\noindent As for the further analysis of the design problem, it would be of interest to derive quantitative
dependence of the optimal observability constant on the positions and stiffnesses of the
connectors, and thereby to determine or approximate the optimal configurations. The existence
result proved here supplies the mathematical basis for such a study, but does not identify a
minimizer explicitly. Another question is to justify the design criterion for connectors of
nonzero width by passing to the limit in the optimal control costs. Finally, to deal with bending
motions or with separation of the two members in the transverse direction, the mechanical model
has to be extended to include the corresponding displacement fields, as in~\cite{MNSS}.
\bigskip

\noindent\textbf{AI Declaration}\par

\noindent AI tools were used to assist with language editing and the optimization of selected proofs. The authors are responsible for all mathematical results, references, and the final version of the manuscript.

\end{document}